\documentclass[12pt, reqno]{article}

\usepackage{hyperref}
\usepackage{theoremref}
\usepackage{amsmath}
\usepackage{amsthm}
\usepackage{amssymb}
\usepackage{amsfonts}
\usepackage{color}
\usepackage[all]{xy} \xyoption{all}
 \usepackage{wrapfig}
\usepackage{enumerate}
\usepackage{enumitem}
\usepackage{slashed}

\usepackage{tikz, tikz-3dplot, pgfplots}
\usetikzlibrary{decorations.markings, patterns, decorations.pathmorphing}

\tikzset{->-/.style={decoration={
markings,
mark=at position #1 with {\arrow{>}}},postaction={decorate}}}

\def\VV{\mathbb{V}}

\def\CC{\mathbb{C}}

\def\EE{\mathbb{E}}

\def\RR{\mathbb{R}}
\def\ZZ{\mathbb{Z}}
\def\QQ{\mathbb{Q}}

\def\gen{\mathfrak{g}}

\def\Fen{\mathfrak{F}}

\def\Len{{\mathfrak L}}

\def\Sen{\mathfrak{S}}

\def\Dc{\mathcal{D}}
\def\Ec{\mathcal{E}}
\def\Fc{\mathcal{F}}
\def\Gc{\mathcal{G}}

\def\Lc{\mathcal{L}}
\def\Mc{\mathcal{M}}

\def\Oc{\mathcal{O}}
\def\Pc{\mathcal{P}}

\def\Sc{\mathcal{S}}
\def\Tc{\mathcal{T}}

 \def\bPhi{{\mathbf{\Phi}}}
\def\bPsi{{\mathbf{\Psi}}}
 \def\bv{{{\bf v}}}

 \def\Aut{{\on{Aut}}}
 
 \def\be{\begin{equation}}
 \def\bef{\begin{figure}}
 \def\bem{\begin{matrix}}
 
 \def\bpm{\begin{pmatrix}}
 \def\Br{{\on{Br}}}
  
 \def\btp{\begin{tikzpicture}}

  \def\Conv{{\on{Conv}}}

 \def\del{{\partial}}

  \def\ee{\end{equation}}
 
 \def\enf{\end{figure}}
 \def\enm{\end{matrix}}
 
 \def\epm{\end{pmatrix}}
 \def\eps{{\varepsilon}}
 \def\etp{\end{tikzpicture}}

 \def\FS{{\on{FS}}}

  \def\Hom{\operatorname{Hom}\nolimits}

 \def\Id{\operatorname{Id}\nolimits}

  \def\Lie{{\on{Lie}}}

   \def\lra{\longrightarrow}

\def\on{\operatorname}
\def\ol{\overline}

\def\oo{{\infty}}

 \def\Perv{{\on{Perv}}}
 \def\phi{{\varphi}}

  \def\Sol{{\on{Sol}}}

 \def\Th{{\on{Th}}}
 
 \def\Tr{{\on{Tr}}}

\def\ul{\underline}

\def\Vect{\on{Vect}}

\def\wh{ \widehat}
\def\wt{\widetilde}

\def\Z2YD{{{1\over 2}\ZZ\hskip -0.07cm\on{YD}}}

\def\0{{\ol{0}}}
\def\1{{\ol{1}}}
\def\(({(\hskip -1mm (}
\def\)){)\hskip -1mm )}
\def\-{{\setminus}}
\def \= {{\,\simeq \,}}
 \def\be{\begin{equation}}
\def\ee{\end{equation}}
\def\ed{\end{document}}

 \numberwithin{equation}{section}
 
\newtheorem{thm}[equation]{Theorem}

\theoremstyle{definition}

\theoremstyle{remark}

\newtheorem{rem}[equation]{Remark}

\newtheorem{exas}[equation]{Examples}

\numberwithin{itemcounter}{subsection}

\begin{document}

\title{  Algebra of the Infrared, secondary polytopes and perverse schobers }

 \author{Mikhail Kapranov\footnote{   Research supported by the  World Premier International Research Center Initiative (WPI Initiative), 
 MEXT, Japan and  by the JSPS  KAKENHI grant 20H01794},
  Yan Soibelman\footnote {Research  supported by  the NSF and the Simons Foundation grants as well
as an IHES Sabbatical Professorship. } }


\maketitle

\begin{abstract}
  This survey paper, based on a talk at the International Congress of Basic Science in
Beijing in July 2025, summarizes
  the joint work of the authors with  M. Kontsevich \cite{KKS} establishing the relation between
  the  ``Algebra of the Infrared" of D. Gaiotto, G. Moore and E. Witten \cite{GMW}
  and the theory of secondary polytopes introduced in the 1990s 
  in the study of higher-dimensional discriminants.

   It also summarizes the subsequent work with   L. Soukhanov
  \cite{KSS}  where the tunnelling data were observed to be similar to linear algebra data describing perverse sheaves on the complex plane except that in the physical context vector spaces are replaced by triangulated categories. The relevant concept here is that of perverse schobers, which are conjectural categorical analogs of perverse sheaves proposed  by M. Kapranov and V. Schechtman \cite{KS}. 
  
  Finally, we sketch a research program of extending these ideas to $4$-dimensional
  theories and the resurgence formalism. 
\end{abstract}

 \tableofcontents

 
 
\numberwithin{equation}{section}
 
  \section {$N=(2,2)$ supersymmetric theories}

The ``infrared limit" analysis of quantum field theories describes them in terms of vacua and of tunnellng effects between the vacua. For 2-dimensional massive (2,2)-supersymmetric theories the set of vacua is typically discrete. In this case D. Gaiotto, G. Moore and E. Witten 
\cite{GMW} introduced a new ``Algebra of the Infrared" which unites various tunnelling phenomena into a single algebraic structure based on the concept of plane webs.

\subsection{Vacua in a theory and the IR limit}
 
 Generally, by
 a {\it vacuum} physicists understand a quantum state (vector) with minimal energy
 i.e., an eigenvector of the Hamiltonian with minimal eigenvalue. Understood as such,
 the ``space of vacua'' is a linear space (the eigenspace with this eigenvalue).
 
 \vskip .2cm
 
 However, there is a different sense in which the expression ``space of vacua''
 is used in modern QFT. It has to do with theories defined on the non-compact
 space-times such as $\RR^{d-1,1}$ or its super-extensions.
 Since the space-time is non-compact, in order to completely define such a theory $\Tc$,
 one needs to specify certain admissible boundary conditions at spatial infinity. 
 Each such  completion then has its own vacuum, a state of minimal energy.

 \vskip .2cm
 
 The  {\em moduli space  of vacua} is a shorthand (synecdoche) term for
 the set $\VV$ of such completions. 
   Endowed with an appropriate topology it becomes a (nonlinear) manifold 
   or maybe a stack. 
 In {\it massive} theories the set of vacua  $\VV$ is 
 discrete, e.g. finite $\VV=\{v_1,\cdots, v_N\}$. 
 
\vskip .2cm
 
The  {\it infrared limit} (with respect to the renormalization group action)
of a theory associates  to it the moduli space of   vacua  as well as certain interaction
 (``tunnelling\rq\rq{}) between them. 
 
 \vskip .2cm
  
In the case of $2d, N=(2,2)$ supersymmetric theories the infrared (IR for short) limit 
 can be described via the 
 {\it Algebra of the Infared} proposed in the seminal paper by Gaiotto-Moore-Witten
 \cite{GMW}. Its aim is to understand  the IR limit  via a certain  algebraic formalism.

\subsection{$N=(2,2)$ supersymmetric theories in 2 dimensions}

  In a 2d QFT on the flat spacetime we work with bosonic coordinates $(x,t)$. 
 The   2d spacetime translations (energy-momentum) $H=i\del_t, P=i\del_x$, 
  part of 2d Poincar\'e algebra, act on everything.
  To allow for $N=(2,2)$ supersymmetry, one adds two fermionic coordinates
  $\xi_1, \xi_2$ and their formal complex conjugates. 
  
  \vskip .2cm

 The fundamental  feature of  the  super-world is that we have square roots of
 spacetime translations. The simplest example is the derivation 
 \[
  \sqrt{\del_x }\,= \,Q\, =\, \del_\xi + \xi \del_x,  \quad   [Q,Q]_+=2\del_x
  \]
on the superspace $\CC^{1|1}$ with cordinates $(x, \xi)$. 

 \vskip .2cm 

In the case of $2d$,  $N=(2,2)$ supersymmetry the corresponding ``roots'' 
form the following {\em supersymmetry algebra}:  
\[
\begin{gathered}{}
 [Q_+, \ol Q_- ]_+ = H+P, \quad [Q_-, \ol Q_+]_+ = H-P 
 \\
 [Q_+, Q_-]_+= c, \quad [\ol Q_+, \ol Q_-]_+ = \ol c.
 \end{gathered}
\]
Here $c\in \CC$ is the {\em  central charge} of the theory, an important quantity. 
Its non-vanishing is considered as an {\em anomaly}, deviation from what is
expected classically.

 \subsection{The $\CC$-plane of central charges}
  We assume that the set of vacua $\VV$ is discrete (this includes all massive theories)
 and, moreover, finite: 
 $\VV=\{v_1,\cdots, v_N\}$. 
 
 For each $i,j=1,\cdots, N$ we then have the 
 {\em  $(i,j)$-sector} of the theory, corresponding to 
  transitioning from $v_i$ to $v_j$. These sectors may have different central charges.
  Denoting by $c_{ij}$ the central charge of the $(i,j)$-sector, we have on
  general grounds \cite{WO} the {\em cocycle condition}
   \[
 c_{ij} + c_{jk} = c_{ik}. 
 \]
 This allows us to write 
  $c_{ij}=w_i-w_j$ with $(w_1,\cdots, w_N)$  defined up to simultaneous shift. 
 Typically, the  $w_i$ can label the
 vacua: $\VV\hookrightarrow \CC$, which we also assume.

Modelling the IR theory categorically,  one assigns to the $i$th vacuum $v_i$  a {\em  local vacuum} (or  {\em ``D-brane''}) {\em category} $\bPhi_i$.   
It is convenient to imagine that  the category $\bPhi_i$ ``sits'' at $w_i\in \CC$. In this language the transport between the vacua $i$ and $j$ is encoded into a 
 {\em transport}  (or {\em tunnelling}) functor $T_{ij}: \bPhi_i \to \bPhi_j$.

\subsection{Example: the Landau-Ginzburg (LG) model}\label{subsec:LG}

 The Landau-Ginzburg (LG) model is a a $2d$ analog of the motion in a potential field. 
The ingredients  are   a  K\"ahler manifold $X$ and a proper holomorphic map 
$W: X\to \CC$  (``superpotential''). The fields are, first of all, a map $\sigma: \RR^{1,1}\to X$
plus various fermionic fields.  This gives a $2d$ $N=(2,2)$ supersymmetric theory. 
The ``effective potential\rq\rq{} 
(obtained by eliminating the fermionic fields)  is equal to $|dW|^2$ and 
 describes the set  $\VV$ of vacua as the set  of critical points of $W$. 
The  LG theory   is massive  if the critical points are non-degenerate (Morse) .

   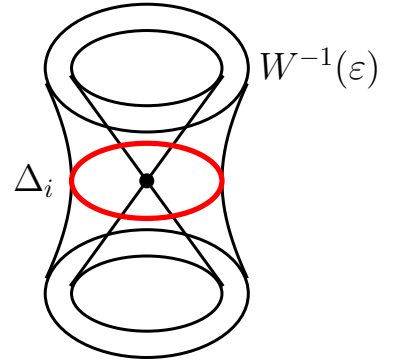
\begin{wrapfigure}{r}{0.25\textwidth}
    \begin{tikzpicture}[scale=0.5]

\node (0) at (0,0){};
\fill (0) circle (0.2);

\draw[ line width=.4mm] (0,3) ellipse (2cm and 1cm);
\draw[ line width=.4mm] (0,3) ellipse (2.7cm and 1.7cm);

 \draw[ line width=.4mm] (0,-3) ellipse (2cm and 1cm);
  \draw[ line width=0.4mm] (0,-3) ellipse (2.7cm and 1.7cm);
\draw [ line width=.4mm] ( -2.035, -2.8) -- ( 2.035, 2.8);
\draw [ line width=.4mm] (2.035, -2.8) -- ( -2.035, 2.8);

\pgfmathsetmacro{\e}{1. 77}   
\pgfmathsetmacro{\a}{2}
\pgfmathsetmacro{\b}{(\a*sqrt((\e)^2-1)}
\draw[line width =0.4mm] plot[domain=-0.8:0.8] ({\a*cosh(\x)},{\b*sinh(\x)});
\draw [line width =0.4mm]  plot[domain=-0.8:0.8] ({-\a*cosh(\x)},{\b*sinh(\x)});

\draw[ color=red, line width=.7mm] (0,0) ellipse (2cm and 1cm);

\node at (4.6,3){\large$W^{-1}(\eps)$};

\node at (-3.,0) {\large$\Delta_i$};

\end{tikzpicture}
\caption{The Lefschetz vanishing cycle.} 
\label{fig:lefschetz}

\end{wrapfigure}

Suppose that $W$ is a holomorphic  Morse function with the set of different critical points $z_1,\cdots, z_N$ and the corresponding set of critical values $w_i=W(z_i)\in\CC$
which we assume distinct. Then the central charge of the $(i,j)$-sector is 
 $c_{ij}=w_i-w_j$.  There are versions of this, including e.g., 
isolated  non-Morse  critical points $z_i$.

In the case of LG theory the $D$-brane category  $\bPhi_i$
can be seen as a categorification of  the vector space  of vanishing cycles
at $z_i$ which we denote $\Phi_i$. More precisely, 
  $\bPhi_i$  is the  {\it  local Fukaya-Seidel category} which categorifies
 $\Phi_i$.
 For example, if $z_i$ is a holomorphic Morse critical point,
then:

 \begin{itemize}
 \item $\Phi_i=\CC\cdot \Delta_i$ is a $1$-dimensional vector space
  spanned by the classical vanishing cycle $\Delta_i$
(sphere) of Lefschetz, see Fig. \ref{fig:lefschetz}, while 

\item $\bPhi_i\= D^b(\Vect_\CC)$ (simplest triangulated category) is generated
by this  sphere (which is a Lagrangian submanifold) as an object of
the local Fukaya category. 
\end{itemize} 

\begin{figure}[h]
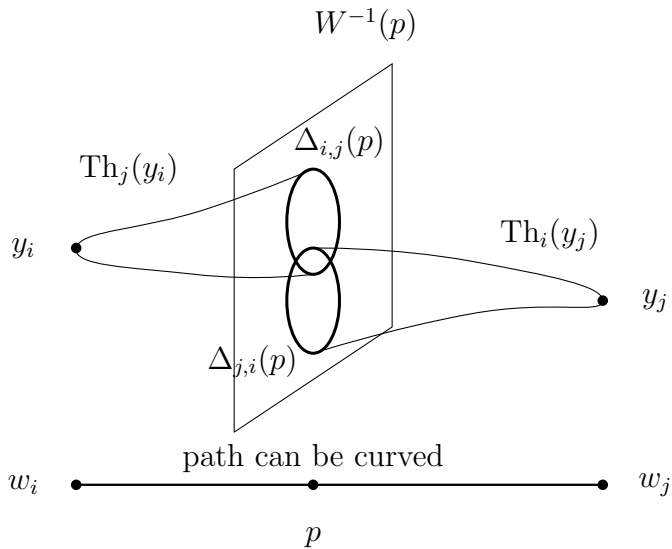


\centering
\btp[scale=0.35]
\node (i) at (-9,0){};
\node (j) at (11,0){};
\node (p) at (0,0){};
\fill (i) circle (0.2);
\fill (j) circle (0.2);
\fill (p) circle (0.2);
\node (xi) at (-9,9){};
\node(xj) at (11,7){};
\fill (xi) circle (0.2);
\fill (xj) circle (0.2);

\draw[ line width=.4mm] (0,7) ellipse (1cm and 2cm);
\draw[ line width=.4mm] (0,10) ellipse (1cm and 2cm);

\draw[line width=0.3mm] plot[smooth, tension=1]  coordinates{ (-9,0)
 (11,0)
};

\draw (-3,2) -- (-3,12) -- (3,16) -- (3,6) -- (-3,2);

\draw plot [smooth, tension=1]  coordinates {(0,12) (-4,10.5) (-9,9) (-4.5, 8) (0,8)};

\draw plot [smooth, tension=1]  coordinates {(0,9) (6,8.5) (11,7) (6,6.5)  (0,5)};

\node at (-11,0) {
$w_i$};
\node at (13,0) {
$w_j$};
\node at (0,-2){
$p$};

\node at (0,1){path can be curved}; 
\node at (-11,9) {
$y_i$};
\node at (13,7) {
$y_j$};

\node at (2,17.5){
$W^{-1}(p)$};

\node at (-7,12){
$\on{Th}_j(y_i)$};
\node at (9,9.5){
$\on{Th}_i(y_j)$};
\node at (1,13){
$\Delta_{i,j}(p)$};
\node at (-2.3,4.7){
$\Delta_{j,i}(p)$};
\etp
\caption{Intersection of Lefschetz thimbles.}
\label{fig:thimbles}
\end{figure}

 \noindent At the level of spaces of vanishing cycles we have   $t_{ij}: \Phi_i \to \Phi_j$ (transport maps). If $W$ is holomorphic Morse function,  then $\Phi_i=\CC\cdot \Delta_i\simeq \CC$. Then 
$t_{ij}: \CC\to \CC$ is   the multiplication with  a number
which is  the intersection number of the corresponding
 {\em Lefschetz thimbles}. That is, we join the points $w_i$ and $w_j$
 by a {\em straight line interval} $[w_i, w_j]$ and form the thimbles
 $\Th_j(y_i)$, $\Th_j(y_j)$ coming out of $y_i$ and $y_j$ respectively
 in the direction of the other point and lying over the interval. Then we restrict
 to the preimage of 
 some middle point $p\in [w_i, w_j]$ and $t_{ij}$ is the intersection index
 of $\Delta_{ij}(p) = \Th_j(y_i)\cap W^{-1}(p)$ and 
  $\Delta_{ji}(p) = \Th_j(y_j)\cap W^{-1}(p)$ inside $W^{-1}(p)$, 
 see Fig. \ref{fig:thimbles}. Note that a priori we can take any curved path in the $w$-plane
 joining $w_i$ and $w_j$ and avoiding other $w_k$, but for the purposes of infrared
 analysis it is important to take the straight interval $[w_i, w_j]$. 
 This has to do with the Fourier transform point of view on tunneling, see below.

 \noindent In the same situation
 $\bPhi_i= D^b(\Vect_\CC)$   and the transport functor  $T_{ij}$ is  the  tensor
 multiplication with the graded vector space $\Hom$ in the Fukaya category 
 of the generic fiber $W^{-1}(p)$:

  \[
 T_{ij}(U)= U\otimes \Hom_{\on{Fukaya}(W^{-1}(p))}(\Delta_{i,j}(a), \Delta_{j,i}(b)).
 \]
 The Euler characteristic of this graded vector space is the intersection index $t_{ij}$,
 so $T_{ij}$ can be seen as a categorification of $t_{ij}$.

 \section{Appearance of secondary polytopes} 

 \subsection{Secondary polytopes in algebraic geometry}

 \begin{wrapfigure}{r}{0.28\textwidth}
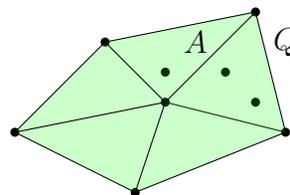

  
\btp[scale=.4, baseline=(current  bounding  box.center)]
\node (1) at (0,0){}; 
\fill (1) circle (0.15);

\node (2) at (4,-1){}; 
\fill (2) circle (0.15);

\node (3) at (-1,-3){};
\fill (3) circle (0.15);

\node (4) at (-5,-1){};

\fill (4) circle (0.15);

\node (5) at (-2,2){};
\fill (5) circle (0.15);

\node (6) at (3,3){};
\fill (6) circle (0.15);

\node (7) at (0,1){};
\fill (7) circle (0.15);

\node (8) at (2,1){};
\fill (8) circle (0.15);

\node (9) at (3,0){}; 
\fill (9) circle (0.15);

\draw (0,0) -- (4,-1);

\draw [draw=black, fill=green, fill opacity=0.2] 
(4,-1) -- (-1,-3) -- (-5, -1) -- (-2,2) -- (3,3) -- cycle;

\draw (0,0) -- (-1,-3); 

\draw (0,0) -- (3,3);
 
\draw (-2,2) -- (0,0); 
 
\draw (0,0) -- (-5,-1);

\node at (1,2){$A$}; 

\node at (4,2){$Q$}; 

\etp

\caption{ A subset  $A\subset \RR^d$, its convex hull $Q$ and a triangulation of $Q$ with vertices in $A$}

\label{fig:secondary}

\end{wrapfigure}

Secondary  polytopes were        introduced by Gelfand-Kapranov-Zelevinsky \cite{GKZ}
 in 1990's. These are convex polytopes $\Sigma(A)$ whose vertices correspond to
 {\em  triangulations} of other, ``primary" polytopes $Q=\on{Conv}(A)\subset \RR^d$,
 see Fig. \ref{fig:secondary}. 
 
 \vskip .2cm

The original motivation came from toric geometry and an attempt  to understand   discriminants of polynomials in $d$ variables. 
That is we consider the set  $A$  of monomials
coded numerically by {\em integer }vectors (e.g. $x^2y^3z^5 \mapsto (2,3,5)$). 
The triangulations   of $Q$ with vertices in $A$ (which are regular in the sense of existence
of convex piecewise-affine functions)  correspond
 some ``extreme''  monomials in the discriminant.  They also correspond to degenerations of
the {\em toric variety} associated to $A$ (used in Mirror Symmetry). 
  In this approach it is important that the points are integer, although the concept
  of the secondary polytope makes sense for any finite set of points $A\subset \RR^d$.

\subsection{Secondary polytopes ``in physics'': factorization and $L_\oo$-algebras }

 In   the paper  \cite{KKS} by M. Kapranov, M. Kontsevich and Y. Soibelman,  a new application of the notion of secondary polytope was proposed.
Here, the set 
  $A=\{w_1,\cdots, w_N\}\subset \CC=\RR^2$ is  the set of critical values of the superpotential $W$.

\vskip .2cm

The important observation is that  faces $F$ of $\Sigma(A)$ (for any $A\subset\RR^d$)
enjoy a natural {\em  factorization property}. That is. 
a face $F=F_\Pc$ corresponds to a polyhedral decomposition $\Pc = \{Q_i\}$, so that 
vertices of such a  face  correspond to further triangulations.  Under genericity assumptions one has:
\[
F_\Pc = \prod \Sigma(Q_i).
\]
Let $C_{\text{polyh.}}(\Sigma(A))$ be the polyhedral chain complex of $\Sigma(A)$ with
complex coefficients. We have a similar complex for any $A'\subset A$. The factorization property
implies that 
 \[
 \bigoplus _{A'\subset A} S^\bullet\bigl(C_{\text{polyh.}}(\Sigma(A'))\bigr)
 \]
  becomes a commutative differential graded algebra (DGA). Further, it is isomorphic to
  the Chevalley-Eilenberg complex
$ C^\bullet_\Lie(\gen_A)$ for
some $L_\oo$, i.e.,  homotopy Lie algebra  $\gen_A$. 
The (higher) 
  Lie brackets in $\gen_A$ are given by  concatenations when possible and are zero otherwise.
  Thus provides an equivalent formulation of the $L_\oo$-algebra of plane webs of 
   Gaiotto-Moore-Witten \cite{GMW} . The webs appear (for $d=2$) as geometric objects
   Poincar\'e dual to triangulations and polygonal subdivisions.

  In fact  the above  definition makes sense in $\RR^d$ for any $d$.

\subsection{$A_\infty$-algebra in the case $d=2$}

 \begin{wrapfigure}{r}{0.2\textwidth}
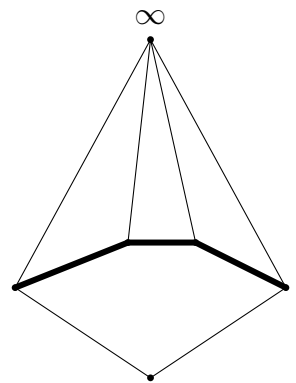

\btp[scale=0.3]
\node (1) at (0,0){};
\node (2) at (6,4){};
\node (3) at (2,6){};
\node (4) at (-1,6){};
\node (5) at (-6,4){};
\node (6) at (0,15){};

\fill(1) circle (0.15); 
\fill(2) circle (0.15); 
\fill(3) circle (0.15); 
\fill(4) circle (0.15); 
\fill(5) circle (0.15); 
\fill(6) circle (0.15); 

\draw (0,0) -- (6,4) -- (2,6) -- (-1,6) -- (-6,4) -- (0,0); 
\draw (-6,4) -- (0,15) -- (-1,6); 
\draw (2,6) -- (0,15) -- (6,4); 

\draw[line width=0.8mm] (-6,4) -- (-1,6) -- (2,6) -- (6,4); 

\node at (0,16){$\oo$ };
\node at (0,3){ };
\node at (-3,6){ };
\node at (0.5,6.7){ };
\node at (4,6){ };
\node at (-2,9){ };
\node at (0.5,9){ };
\node at (2.4,9){ };

\node (d) at (10,7){ };
\node (e) at (6,6){};

\node (f) at (8,1){ };
\node (g) at (4,2){};

\etp
\caption{ The extended polytope $\Conv(A\cup \{\oo\})$.  }
\label{fig:pos-neg-bou}
\end{wrapfigure}

 In the case of $\RR^2$ 
we can add one point
at infinity, getting an ``extended polytope'' $\wt Q= \Conv(A\cup \{\oo\})$
with one vertex at $\oo$., see
See Fig.\ref {fig:pos-neg-bou}. The infinite triangles of any triangulation of $\wt Q$
are then ordered anticlockwise. This natural ordering allows us to define, similarly
to the above,   not just homotopy commutative but  a
{\rm homotopy associative}, or $A_\oo$-algebra
$R:=R_{A\cup \{\oo\}}$  of infinite polygons. Further, possible decompositions of infinite
polygons into finite and infinite ones define
 a morphism of $L_\oo$-algebras  from
$\gen_A$   to $\on{Hoch}^\bullet(R_{A\cup\{\oo\}})$
(the deformation complex). The main result of \cite{KKS} is:

\begin{thm}\label{thm:qis}
This  morphism is a
quasi-isomorphism.

\end{thm}

 The general deformation theory says that in order to 
define a deformation of $R$ we need to specify the corresponding
Maurer-Cartan element in  $\gen_A$. 
From this perspective, in the work of Gaiotto-Moore-Witten some particular MC element was proposed, and it was conjectured that such a choice gives the
Fukaya-Seidel category.

\section{Perverse sheaves and perverse schobers}

\subsection{Analogy with perverse sheaves}\label{subsec:an-perv}

 Let $(X,\Sc)$ be a  Whitney stratified complex manifold. We then have the abelian
category $\Perv(X, \Sc)$ formed by perverse sheaves (of $\CC$-vector spaces,
with middle perversity)
on $X$  constructible  w.r.t. $\Sc$. They are defined among all $\Sc$-constructible
complexes by 
 certain dimensional conditions on supp$(\ul H^i(\Fc^\bullet))$. 

By the Riemann-Hilbert correspondence in the form proposed by Kashiwara and Mebkhout in the 1980\rq{}s there is an equivalence of the category of perverse sheaves and the category 
of holonomic regular {$\Dc$-modules}   $\Mc$, via
$\Fc^\bullet = {R\ul \Hom}_\Dc(\Mc, \Oc)$,  the derived sheaf of solutions. 

There are several cases when the categories $\Perv(X,\Sc)$
can be identified explicitly, as category of representations of quivers with relations.

\begin{exas}\label{exas:Perv}
 (a)  On $(\CC,0)$  (i.e., $X=\CC$ and $\Sc$ consists of $0$ and $C\- \{0\}$). 
One can also consider the case when $X=U$ is a small disk around a point $p$
and $\Sc$ consists of $p$ and $U\-\{p\}$. In any such case
$\Perv(\CC,0)$ is identified with the category of diagrams of vector spaces
 $  \xymatrix{ 
 \Phi \ar@<.4ex>[r]^a&\Psi \ar@<.4ex>[l]^{a'}}$ such that  $1_\Psi - aa'$ and $1_\Phi-a'a$ are invertible (these are the monodromies). 
 Here $\Phi$ is called the space of {\em vanishing cycles} and    $\Psi$
 the space of {\em nearby cycles}, or the {\em generic fiber}. The space $\Phi_i$ of
 vanishing cycles in \S \ref{subsec:LG} is a particular case of this construction
  to a certain perverse sheaf
 ({\em Lefschetz perverse sheaf}) on $\CC$, associated to $W$ and considered
 in a neighborhood  $U$ of $w_i$, see \S \ref{subsec:schob-LG}
 below.

 \  \vskip .1cm
 
 (b) On $(\CC, w_1,...,w_N)$  (i.e., $X=\CC$ and $\Sc$ consists of the $w_i$ and
 $\CC\-\{w_1,\cdots, w_N\}$). Here the category $\Perv$ can be described
 by ``amalgamation'' of the descriptions near each $w_i$ given by (a),
 as shown by 
  Gelfand-MacPherson-Vilonen \cite{GMV}. More precisely, $\Perv(\CC, w_1,
 \cdots, w_N)$ is identified with the category of diagrams of vector spaces
   \[
\xymatrix{\Phi_N
\ar@<.4ex>[rd]^{a_N}
&
\\
\vdots & \Psi \ar@<.4ex>[ul]^{a'_N}
\ar@<.4ex>[dl]^{a'_1}
\\
\Phi_1 \ar@<.4ex>
[ur]^{a_1}
}
\label{eq:GMV-sh}
\]
where each  
$  \xymatrix{ \{ \Phi_i \ar@<.4ex>[r]^{a_i}&\Psi \ar@<.4ex>[l]^{a'_i}\}}$
is as in (a). 

 \vskip .1cm
 
  (c) Further, in the context of (b) one can describe the category
 \[
 \ol\Perv(\CC, w_1,\cdots, w_n) \,=\,\Perv(\CC, w_1,\cdots, w_n)/\on{Const},
 \]
 the quotient by the Serre subcategory of constant sheaves. Here only the $\Phi$-data
 survive and the result of \cite{GMV} identifies this category with the
 category of ``matrix diagrams'' consisting of vector spaces $\Phi_1, \cdots, \Phi_n$,
 {\em transport maps} $t_{ij}: \Phi_i\to \Phi_j$, $i\neq j$ and monodromy
 automorphisms $\mu_i\in \Aut(\Phi_i)$. 
 The transport maps appear as $t_{ij}= a_j' a_i$.
 \end{exas}
 
  Note that one can define the  transport maps $t_{ij}( \alpha)$ along any path 
  $\alpha$ from $w_i$ to $w_j$,
 avoiding other $w_k$. It is the composition of the map $a_i: \Phi_i\to\Psi_i$,
 (where we think of $\Psi_i$ as the stalk at a point of $\alpha$ near $w_i$,
 then the monodromy $\Psi_i\to\Psi_j$ along $\alpha$ and then $a'_j:\Psi_j\to\Phi_j$,
 see Fig. \ref {fig:transport}. This generalizes the construction of 
  \S\ref{subsec:LG}. The dependence on $\alpha$ is governed by the
 {\em abstract Picard-Lefschetz identities} \cite{KSS}. 
 The definition $t_{ij}= a_j' a_i$ in the equivalence of  \cite{GMV} corresponds to some
 particular path, see \S \ref {subsec:schob-spid}  below. 
 In 
 \S\ref{subsec:fourier} we will use {\em  rectilinear
  transport maps} $t^R_{ij} = t_{ij}[w_i, w_j]$ taken along
 the straignt interval $[w_i, w_j]$ (assuming there are no other $w_k$ on this interval).

 \bef[h]
\centering
\btp [scale=.8, baseline=(current  bounding  box.center)]

\node (i) at (0,0){};
\fill (i) circle (0.1);

\node (j) at (12,0){};
\fill (j) circle (0.1);

\node (k) at (5,-2){};
\fill (k) circle (0.1);

\node (b) at (4,0.3){$$};
\fill (b) circle (0.08);

\draw [dotted, line width =0.3mm]  (0,0) circle (1cm);
\draw [dotted, line width =0.3mm]  (12,0) circle (1cm);

\draw[->, line width = 0.2mm] plot [smooth,tension=1.5] coordinates{
(0,0) (4,0.3)
(8,-0.3)
(12,0)
};
\node at (-0.5,0){$w_i$};
\node at (12.5, 0){$w_j$};
\node at (5.5, -2){$w_k$};
\node at (9,-0.7){\large$\alpha$};
\node at (4,-0.3){$p$};
\node at (4,0.8){$\Psi_\alpha$};
\node at (1.4, 0.7){$\Phi_{i,\alpha}$};
\node at (10.5, 0.3){$\Phi_{j,\alpha}$};

\etp

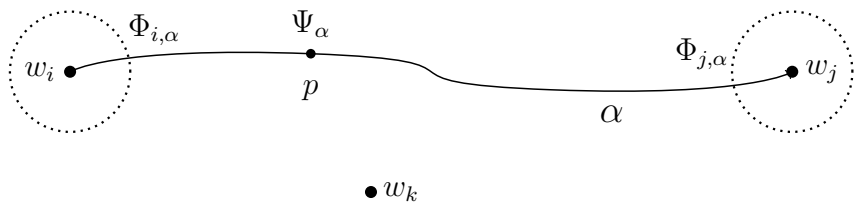
\caption{ The transport map.  }
\label{fig:transport}
\enf

\vskip .2cm

   The observation of Kapranov-Soibelman-Soukhanov \cite{KSS} is that this is very  similar to
 the tunnelling data
 of the infrared formalism. So one can formulate  the main thesis of \cite{KSS}
  as the following slogan: 
 \[
 \boxed{
 \begin{gathered}
 \text{GMW (Gaiotto-Moore-Witten formalism)} = 
 \\
 =\text{GMV (Gelfand-MacPherson-Vilonen
 classification)}. 
 \end{gathered}
 }
 \]

\subsection{Perverse schobers}

 Kapranov and Schechtman \cite{KS} proposed in 2014 a categorification of perverse sheaves, which they called {\it perverse schobers}.
 The data consisting of vector spaces is replaced
  in their approach   by the  data consisting of triangulated categories.
 The general definition of perverse schobers is still unclear.  So far people gave the definition in some special cases by directly categofying various quiver
 descriptions. 
 
 \begin{exas}
 
 (a)  On $(\CC,0)$  a perverse schober is defined by the data of 
a pair of  spherical functors
$  \xymatrix{ 
 \bPhi \ar@<.4ex>[r]^{a}&\bPsi \ar@<.4ex>[l]^{a^*}
 }$, 
 i.e., an  adjoint pair s.t.
 $\on{Cone}\{aa^*\to \Id_\bPsi\}$, $\on{Cone}\{\Id_\bPhi\to a^*a \}$  are equivalences.
 This naturally categorifies  Example \ref{exas:Perv}(a).

  \vskip .1cm
 
  (b) More generally on $(\CC,w_1,...,w_N)$ we have similar data which are categorical analogs of  Example   \ref{exas:Perv}(b), namely a datum of
 triangulated categories and adjoint pairs of functors 
 \be\label{eq:GMV-schob}
\xymatrix{\bPhi_N
\ar@<.4ex>[rd]^{a_N}
&
\\
\vdots & \bPsi \ar@<.4ex>[ul]^{a^*_N}
\ar@<.4ex>[dl]^{a^*_1}
\\
\bPhi_1 \ar@<.4ex>
[ur]^{a_1}
}
\ee
 where  each $a_i$ is  spherical. This gives
{\em transport functors}  $T_{ij}= a_j^* a_i: \bPhi_i\to\bPhi_j$. 
 \end{exas}

 \subsection{A dictionary: theories and schobers}
 
  In \cite{KSS} we proposed the following dictionary. 
 
    \vskip .2cm

    \centerline{ Theories  ($2$d, $N=(2,2)$ etc.) { $\buildrel \text{IR}\over\rightarrow$} Schobers on $\CC= \{\text{central charges}\}.$}
   
   \vskip .2cm
   
   \centerline{ Vacua of a theory   { $\leftrightarrow$} Singular points of the  schober . }
   
      \vskip .2cm
    
  \centerline{  Local D-brane categories  {$\leftrightarrow$} $\bPhi_i$ (vanishing cycle categories). }
  
     \vskip .2cm
    
  \centerline{  Tunneling  between  vacua  {$\leftrightarrow$}  
  Transport functors $T_{ij}$ given by schobers.  }
  
     \vskip .2cm
  
   \centerline{   Infrared Algebra of GMW  { $\leftrightarrow$} {Fourier transform for schobers. }}
   
     \vskip .2cm
     
   \noindent 
   More precisely (see below), the Infrared Algebra
   corresponds to the 
   explicit description of the analog of Stokes data for the Fourier transform. 
     This agrees with the following thesis due to G. Moore: 
      {\em The Infrared Algebra  provides  a  categorification of  the quasi-classical expansion   for}
    \be\label{eq:exp-int}
    \int e^{{i\over \hbar} W(x)} dx. 
    \ee

    \subsection{Perverse sheaves and schobers in Landau-Ginzburg models}
    \label{subsec:schob-LG}
    
 A proper holomorphic map          
$W: X\to\CC$ as before  gives rise to a natural  complex of sheaves $RW_* \ul \QQ_X$ on $\CC$,
the  direct image of constant sheaf $\ul \QQ_X$. It describes the behavior of possible integration
cycles $\Gamma$, and hence the monodromy of the 
 ``fiberwise integral'', i.e., of the 
    mulltivalued  function  $L_W$  on $\CC$ defined by  
    \[
     a \,\,\mapsto \,\, L_W(a) = \int_{\Gamma'\subset W^{-1}(a)} {dx\over dw}, \quad 
     \quad w=\text{coordinate in } \CC. 
     \]

    The function $L_W$ can be seen as an object of categorical level $0$. At the categorical level
   $1$ we have the complex $RW_* \ul \QQ_X$ or, which is better, the perverse sheaf
   $\Lc_W$ defined as the perverse cohomology (truncation using the perverse t-structure) 
      \[
    \Lc_W= \ul H^0_{\on{perv}} RW_* \ul \QQ_X \in \Perv(\CC, w_1,...,w_N).
    \]
 We call $\Lc_W$ the 
  {\em  Lefschetz perverse sheaf}. Its GMV data are as follows. The $\Phi_i$ are
  the classical  vanishing cycles spaces 
    $\Phi_i$  of Lefschetz, see \S \ref{subsec:LG},  while  $\Psi=H^{\on{mid}}(W^{-1}(\text{generic } a\in\CC))$ is the middle cohomology of the generic fiber. Similarly, the transport maps
    $t_{ij}$  for $\Lc_W$ are the 
   classical transport maps  $t_{ij}$ as defined in \S \ref{subsec:LG}, e.g., $t_{ij}$ is
   the intersection number of Lefschetz thimbles, if $W$ is Morse. 
       
       
        Next, at the categorical level 2, we have a natural 
      {\em Lefschetz perverse schober} $\Len_W$  defined as follows. Its categories
      $\bPhi_i$ are the  local Fukaya-Seidel S categories,  while $\bPsi$ 
      is the  Fukaya category of the generic fiber $W^{-1}(a)$. The       
   transport functors $T_{ij}$ are defined  via $\Hom$-complexes in the Fukaya category,
    as before.

\subsection {Fourier transform for $\Dc$-modules and perverse sheaves}
\label{subsec:fourier}

      Note that the exponential integral
     \eqref{eq:exp-int} is nothing but  the $1$-dimensional Fourier transform
     of the fiberwise integral  $L_W$, if we think of this Fourier transform as a function of  $z=1/\hbar$. This suggests  that Fourier transforms (in the appropriate sense)
      of the 
     perverse sheaf $\Lc_W$ and the schober $\Len_W$ are important objects. 
     We start at the level of perverse sheaves. 
     
     \vskip .2cm

Let $D=\CC\langle w, \del_w\rangle$  be the algebra of polynomial differential
operators on $\CC$.   By the Riemann-Hilbert
correspondence,  the category if   holonomic regular left $D$-modules  $M$ 
is equivalent to that of perverse sheaves $\Fc$ on $\CC$ via the functor 
$M\mapsto \Fc=  \Sol(M)= \ul{R\Hom}_D(M, \Oc_\CC)$. 
So starting from $\Fc\in \Perv(\CC, A)$,  $A=\{ w_1,\cdots, w_N\}$,  we uniquely reconstruct $M$.

Now to $M$ we associate its {\em formal Fourier transform} 
  $\wh M $  which is the same $M$ but  considered as a module  over
   $\CC\langle z, \del_z\rangle$,  with
   $z$ acting as $\del_w,$ and $\del_z$ as  $-w$. As a  $\CC\langle z, \del_z\rangle$-module,
   $\wh M$ is still holonomic but typically not regular and has only  two singularities: 
    at $0$, regular and at $\oo$,
   typically irregular.    
   The complex $\wh\Fc = \Sol(\wh M)$
   is still perverse, so  $\wh\Fc\in \Perv(\CC,0)$.
   So $\wh\Fc|_{\CC\- \{0\}} = \Ec[1]$ where $\Ec$ is a local system on $\CC-\{0\}$
   (whose sections are exponential integrals of multivalued solutions of $M$).

 \begin{wrapfigure}{r}{0.4\textwidth}
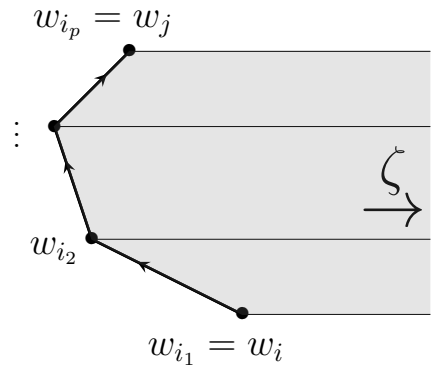
 
 \btp[scale=0.5]
\node at (-4,3){$\bullet$};
\node at (-6,1){$\bullet$};
\node at (-5, -2){$\bullet$};
\node at (-1, -4){$\bullet$};
\draw (-4,3) -- (4,3);
\draw (-6,1) -- (4,1);
\draw (-5, -2) -- (4, -2);
\draw (-1,-4) -- (4, -4);

\draw [line width= 1] (-4,3) -- (-6,1) -- (-5,-2) -- (-1, -4);

\draw[decoration={markings,mark=at position 0.7 with
{\arrow[scale=1,>=stealth]{>}}},postaction={decorate},
line width=1]  (-6,1) -- (-4,3); 

\draw[decoration={markings,mark=at position 0.7 with
{\arrow[scale=1,>=stealth]{>}}},postaction={decorate},
line width=1] (-5, -2) -- (-6,1); 

\draw[decoration={markings,mark=at position 0.7 with
{\arrow[scale=1,>=stealth]{>}}},postaction={decorate},
line width=1] (-1,-4) -- (-5, -2);

\node at((-4.7, 3.7){\large$w_{i_p}=w_j$};
\node at((-7, 1){$\vdots$};
\node at((-6, -2.3){\large $w_{i_2}$};
\node at((-1.7, -5){\large$w_{i_1}=w_i$};

\node at (3, -0.5){\huge$\buildrel \zeta\over\to $};
\filldraw [color=gray, opacity=0.2] (4,3) -- (-4,3) -- (-6,1) -- (-5, -2) -- (-1, -4)-- (4, -4);
\etp
\caption{A $\zeta$-convex polygonal path.} 
\label{fig:path}
\end{wrapfigure}

 More precisely,  $\wh \Fc$ 
    is described
   as in  Example \ref{exas:Perv}(a) by a diagram
   $  \xymatrix{ 
\wh \Phi \ar@<.4ex>[r]^{\wh a}&\wh \Psi \ar@<.4ex>[l]^{\wh a'}}$ whose spaces
are identified as follows \cite{malgrange}: 
$
\wh\Phi$ is identified with  $\Fc_a[-1]$, i.e.,   the generic stalk of the local system  $\Fc[-1]$
on $\CC\- A$, while 
 $\wh \Psi$, i.e., the generic stalk of the local system $\Ec$, 
 is identified with  $\bigoplus\Phi_i$, the direct sum of all the spaces of vanishing cycles
 for $\Fc$.  The identification is not canonical.

Further, $\Mc$ being  irregular, $\Ec$ carries the {\em Stokes structure}
\cite{malgrange, deligne-docs}. It was
  studied and described in various forms in  \cite{malgrange, mochizuki, dagnolo-sabbah}
    and other papers.    In \cite{KSS} we proposed 
    a different form of answer, resembling the Algebra of the Infrared.

  More precisely, the Stokes structure of the kind appearing in this problem
 is determined by one {\em Stokes matrix} $C: \bigoplus \Phi_i \to\bigoplus \Phi_i$,
 see  \cite{deligne-var, KSS}.  Here we fix a direction $\zeta$ at $\oo$ (thought of
 as a complex number, $|\zeta|=1$) and consider the numbering of 
 $A=\{w_1,\cdots, w_N\}$ so that the order increases in the direction
 orthogonal to $\zeta$. Then  $C$ is  block upper-triangular with identities on the diagonal,
 and we need to find the matrix elements $C_{ij}: \Phi_i\to\Phi_j$, $i<j$. 
 
 For this we consider  {\em $\zeta$-convex polygonal paths} out of elements of $A$,
 i.e., sequences
 $w_i=w_{i_1}, w_2, \cdots, w_{i_p}=w_j$ so that each $w_{i_\nu}$ is a vertex
 of the semi-infinite polygon $\Conv(\bigcup  (w_{i_\nu} +\zeta\cdot \RR_+))$,
 see Fig. \ref{fig:path}.   
 
  \begin{thm}\label {thm:baby}\cite{KSS} Suppose no three points $w_i, w_j, w_k$ lie
 on a line in $\CC=\RR^2$. Then  
 \[
 C_{ij} = \sum_{\text{such paths}}  t^R_{i_{p-1}i_p}\circ \cdots \circ  t^R_{i_1i_2}
 \]
 where $t^R_{kl}$  is the  rectilinear transport map from $\Phi_k$ to $\Phi_l$,
 see \S \ref {subsec:an-perv}. 
 \end{thm} 
 
 \noindent  An expression  of this type  can be called as the {\em Baby} (or {\em decategorified}) 
 {\em  Algebra of the Infrared}. So {\em the Baby AIR describes 
 the Stokes matrix of the Fourier transform}.
 
  \begin{rem}
 The appearance of rectilinear transports and convex geometry here  becomes
 natural if we recall that the Fourier transform in the complex domain
 is governed by the domains of decay of the integrands, and that
 exponential functions decay in half-planes. 
 \end{rem}

 \subsection {Schobers and spiders}\label{subsec:schob-spid}

\begin{wrapfigure}{r}{0.3\textwidth}
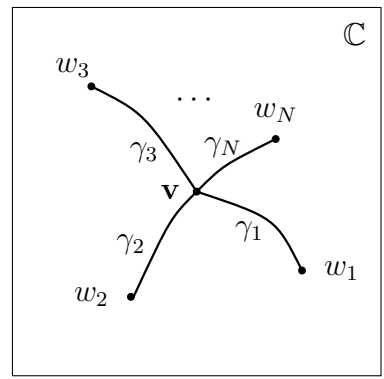
 
\btp[scale=.35, baseline=(current  bounding  box.center)]

\node (1)  at (4, -3){};
\fill (1) circle (0.15);
\node at (5.5, -3) {$w_1$};
\draw[line width = .3mm] (0,0) .. controls (3, -1) .. (4, -3);
\node at (2, -1.5){$\gamma_1$};

\node (2)  at (-2.5,-4){};
\fill (2) circle (0.15);
\node at (-4,-4){$w_2$};
\draw[line width = .3mm] (0,0) .. controls (-1,-1) ..  (-2.4,-4);
\node at (-2.5,-2){$\gamma_2$};

\node (3) at (-4,+4){};
\fill (3) circle (0.15);
\node at (-4.7, +4.7){$w_3$};
\draw[line width = .3mm] (0,0) .. controls (-2,+3) ..  (-4,+4);
\node at (-2,+1.5){$\gamma_3$};

\node (n)  at (3,+2){};
\fill (n) circle (0.15);
\node at (3,+3){$w_N$};
\draw[line width = .3mm] (0,0) .. controls (1,+1) .. (3,+2);
\node at (1,+1.7){$\gamma_N$};

\node (v) at (0,0){};
\fill (v) circle (0.15);
\node at (-1,0){$\bv$};

\node at (0,+3.5){$\cdots$};

\node at (6,6) {$\CC$};
\draw (-7,-7) -- (-7,7) -- (7,7) -- (7,-7) -- (-7,-7);
\etp
\caption{A spider}
\label{fig:spider}

 \end{wrapfigure}
 
  The GMV description \cite{GMV}  of $\Perv(\CC, A)$, $A=\{ w_1,\cdots, w_N\}$ as well as the
 definition of schobers with singularities in $A$ is based  on the choice of a
 {\em spider}, a system $\Sen$ of non-intersecting
 simple curves joining each $w_i$ with a fixed
 point $\bv$. This is a standard concept in singularity theory.
 The transport map $t_{ij}=a'_j a_i$ is, in the notation
 of \S \ref{subsec:an-perv} $t_{ij}(\alpha)$ where $\alpha$ is the path from $w_i$ to $w_j$
 through $v$ in the spider. 
 Note that such  path may be roundabout and not the most convenient: for instance,
 the straight interval $[w_i, w_j]$ is physically more natural. The Algebra of the
 Infrared can be seen as  describing the interaction of these two choices.

 The set of isotopy classes of spiders (with fixed $A, \bv$)  is a torsor over the braid group $\Br_N$ and \cite{GMV} gives explicit formulas (``mutations'') for the change of
 the GMV-data under the  action of the generators of $\Br_N$. These formulas
 can be deduced from the abstract Picard-Lefschetz identities \cite{KSS}. 
 
 These constructions generalize to schobers. 
 The definition is given via a fixed spider
 and can be re-calculated with respect to any other spider using the categorified
 mutations which, in their turn, are defined using {\em abstract Picard-Lefschetz
 triangles}, categorifying the Picard- Lefschetz identities, see \cite{KSS}. 
 
  If a schober $\Fen$ is given in terms of a given spider $\Sen$ by a diagram
 \eqref {eq:GMV-schob}, then the
  transport functors  $T_{ij}= a_j^* a_i$ form  a {\em monad}, or an ``algebra of functors'':
  we have associative compositions
  \[
  T_{jk} T_{ij}  =  a_k^* a_j a_j^* a_i \lra a_k^* a_i = T_{ik}
   \]   
   coming from the counit of the adjunction
    $a_j a_j^*\to\Id$.

      \subsection { The Fukaya-Seidel category with coefficient in a schober $\Fen$}
      
         \begin{wrapfigure}{R}{0.35\textwidth}
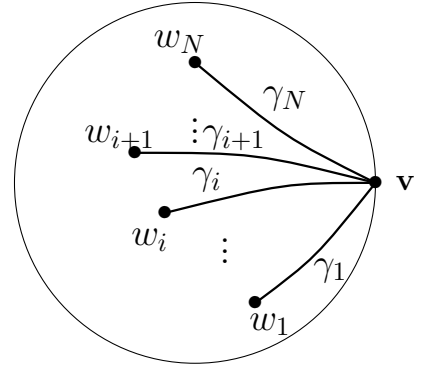
 
   \btp[scale=.4, baseline=(current  bounding  box.center)]
\draw (0,0) circle (6cm);

\node at (6,0){$\bullet$};
\node at (2,-4){$\bullet$};
\node at (-1,-1){$\bullet$};
\node at (-2,+1){$\bullet$};
\node at (0, +4){$\bullet$};
\node at (7,0){$\bv$};

\draw[line width = .3mm] (6,0) .. controls (4, -2.5) ..  (2,-4);
  \draw[line width = .3mm] (6,0) .. controls (3,0) ..  (-1,-1);
   \draw[line width = .3mm] (6,0) .. controls (2,+1) ..  (-2,+1);
 \draw[line width = .3mm] (6,0) .. controls (3,+1.5) ..  (0,+4);

\node at (2.5, -4.7){\large $w_1$};
\node at (4.5, -2.9){\large$\gamma_1$};
\node at (1, -2){\large$\vdots$};
\node at (-1.5, -1.9){\large$w_i$};
\node at (0.4, +0.1){\large$\gamma_i$};
\node at (-2.5, +1.5){\large$w_{i+1}$};
\node at (0, +2){\large$\vdots$};
\node at (1.3, +1.5){\large$\gamma_{i+1}$};
\node at (-0.5, +4.7){\large$w_N$};
\node at (3, +2.8){\large$\gamma_N$};
\etp
\caption{The infinite spider}
\label{fig:inf-spid} 
\end{wrapfigure}

   This construction is analogous to that of D-brane categories in physics.
  Like these, it   
   depends  on the  choice of a half-plane or, what is equivalent, a direction $\zeta\in S^1_\oo$
   near $\oo$. Or, equivalently, we choose a  point
   $\bv$ in   $\CC$ which is very far away (``Vladivostok'') or,
   if we work with a disk $D$ instead of $\CC$, we choose       $\bv\in \del D$. 
   Then, we construct the {\em Infinite  spider} $\Sen_\zeta=\Sen_\bv$
   out of paths going from the $w_i$ straight to $\bv$ (we assume that these
   paths do not intersect), see Fig. \ref{fig:inf-spid}. 
   
    This spider gives rise to the {\em Fukaya-Seidel monad }
   $T(\zeta) = (T_{ij}(\zeta))_{i<_\zeta j}$ (we disregard the remaning $T_{ij}(\zeta)$). 
The {\em Fukaya-Seidel category}   $\FS(\Fen,\zeta)$ is defined as the 
category  of algebras  over the monad $T(\zeta)$. 
This is a  triangulated category with semi-orthogonal
 decomposition $\langle  \bPhi_1, ..., \bPhi_N\rangle$, with
 $w_i$ numbered according to $\leq_\zeta$. 
 
 By construction,  $\FS(\Fen,\zeta)$  is a direct categorification of
 the  Fourier-Sato transform of a perverse sheaf $\Fc$ which
 involves local cohomology $R\Gamma_{\text{half-plane}}(\Fc)$. 
 In other words, $\FS(\Fen,\zeta)$  should be thought of as
 the fiber of the ``schober Fourier transform''  $\wh\Fen$ at the point $\zeta\in S^1_\oo$. 
 
 \subsection {Infrared complex and the Fukaya-Seidel monad}
 
    Let  $\Fen $  be a schober on $\CC$ with singularities in $A=\{w_1,\cdots, w_M\}$.
 We assume that no three points of $A$
   lie on a line. 
   Using Picard-Lefschetz triangles,
   we can define transport functors $T_{ij}(\alpha): \bPhi_i\to\bPhi_j$
   for any path $\alpha$ joining $w_i$ and $w_j$, avoiding other $w_k$,
   as in \S \ref {subsec:an-perv}. In particular, we have rectilinear transport
   functors $T_{ij}^R$ corresponding to $\alpha = [w_i, w_j]$.
   
   Let    $\zeta\in S^1_\oo$ a direction at $\oo$. As in Theorem \ref {thm:baby},
   we define
    iterated rectilinear transport functors $T_\gamma$
    for $\zeta$-convex polygonal  paths $\gamma$.
    They  also  form a monad $R = R(\zeta)$ with 
   $T_\gamma T_{\gamma'}\buildrel =\over\to T_{\gamma\circ\gamma'}$ 
   given by  concatenation, if 
   $\gamma\circ\gamma'$ is $\zeta$-convex
    and by $0$ if  not.       
  This monad is the schober  analog of the $A_\oo$-algebra $R_\oo$ of \cite{GMV}. 
  Further, the following result is
   a schober version of the IR algebra approach of \cite{GMW}.
   It is directly iformualted  in terms of  a categorical analog of the Stokes filtration on the local system of Fukaya-Seidel categories
    $\FS(\Fen, \zeta)$. 
  
   \begin{thm} \cite {KSS} 
   (1)  There exists a differential $d$ in $R(\zeta)$ s.t. $d^2=0$ in the derived
   category.
   
   \vskip .1cm
     
     (2) There further exists a Postnikov system (``filtration'') for $T(\zeta)$ whose 
     ``associated graded object'' 
     is $R(\zeta)$ and the ``connecting morphisms" (in derived category) are
     given by $d$. 
     \end{thm}
     
     The differential $d$ can be called the {\em infrared differential}.
     It is constructed from arrows in various
     Picard-Lefschetz triangles. 
     It should lift to a Maurer-Cartan element in $R(\zeta)$ at dg-level
   (modulo more foundational work on schobers). 
   
 \subsection {Further aspects: $L_\oo$-algebra and secondary polytopes}
  Assuming   $A=\{w_1,..., w_N\}$  sufficiently generic (including no 3 points on a line)
 we construct the secondary polytope $\Sigma(A)$. 
 
 For any subset $A'\subset A$ with $Q'=\Conv(A')$    we define the dg-vector space
 $\EE(A') = \Tr (\text{transport along } Q')$.  Here for an endofunctor $F$
 its categorical trace $\Tr(F)$ is defined as the Hochschild homology of $F$,
 i.e., the dg-space of natural transformations from $\Id$ to $F$.

 \begin{thm} \cite{KSS}
 The sum $\gen_A$ of such $\EE_{A'}$ over all $A'\subset A$
  has a structure of a $L_\oo$
  algebra and we have a quasi-isomorphism  from $\gen_A $ to
  the ordered Hochschild (deformation) complex
  $ C^\bullet (R(\zeta),R(\zeta))$.  
 \end{thm}
 
  So any deformation of $R(\zeta)$ (such as $T(\zeta)$ giving the FS category) comes from a Maurer-Cartan  element
 in $\gen_A$. 
 This is a schober generalization of Theorem \ref{thm:qis}. 
 
 \section{Work in progress: 4d theories and resurgence}
 
  The above approach can be modified to apply  to $N=2$ theories in 4 dimensions
 \cite{GMN-4dWC, GMN-WKB, GMN-2d-4d}.  More precisely, we need to modify
 the approach to apply 
   to resurgence-like phenomena. Indeed, 
 cluster transformations
 (coordinate changes) from \cite{Kont-S} appear as Stokes data for some non-linear
 differential equations (formulated in terms of a Riemann-Hilbert problem). Now, we
 can pass from a non-linear coordinate change to a linear transformation
 on the ring of regular functions (or its appropriate completion).  This means that
 we need to consider  ``irregular linear differential equations of infinite rank'' equipped with
  commutative  algebra structures  on the equations themselves, 
 in particular, in the fibers of the corresponding local systems of solutions. 
 Then the ``Stokes matrices'' will be algebra homomorphisms, i.e., coordinate
 changes. 
 
 \vskip .2cm
 
\begin{figure}[h]
 \centering
 \begin{tikzpicture}[scale=0.6]
 \draw (-4,-4) -- (-4,4) -- (4,4) -- (4, -4) -- (-4, -4); 
  \draw (3,4) -- (3,3) -- (4,3); 
 
 \node at (3.5,3.5) {\large$w$}; 
 
  \node at (1,2 ){\small$\bullet$}; 
   \node at (3,0.3){\small$\bullet$}; 
      \node at (-2,1.5){\small$\bullet$}; 
      \node at (-1,2.5){\small$\bullet$};    
    \node at (-2,-2){\small$\bullet$};    
      \node at (3,-1.5){$\small\bullet$};     
\node at   (0, -2.5){$\cdots$};  
\node at (0,0){(Singularities)} ;  
 \end{tikzpicture}
 \begin{tikzpicture}[scale=0.6]
 \draw[color={white}]    (-4,-4) -- (-4,4) -- (4,4) -- (4, -4) -- (-4, -4);
 \draw   [ decoration={markings,mark=at position 0.99 with
{\arrow[scale=1.5,>=stealth]{>}}, 
mark=at position -1with
{\arrow[scale=1.5, >=stealth]{<}}
},postaction={decorate},
line width=.4mm] (-3.5,0) -- (3.5,0);
\node at (0, 0.5){Fourier}; 
\node at (0, -0.5){transform}; 
 \end{tikzpicture}
 \begin{tikzpicture}[scale=0.6]
 \draw (-4,-4) -- (-4,4) -- (4,4) -- (4, -4) -- (-4, -4);
 \draw (3,4) -- (3,3) -- (4,3); 
  \node at (3.5,3.5) {\large$z$}; 
  \node at (0,0){$\bullet$}; 
 \node at (-0.5, -0.5){$0$}; 
 \filldraw[opacity = 0.2] (0,0) -- (-2,4) -- (2, 4) -- (0,0); 
 \draw (-2,4) -- (0,0) -- (2,4); 
 \node at (0, -1.5){(Sectors)}; 
 \end{tikzpicture}
 \caption{The Borel/Fourier transform}
 \label{fig:borel}
 \end{figure}
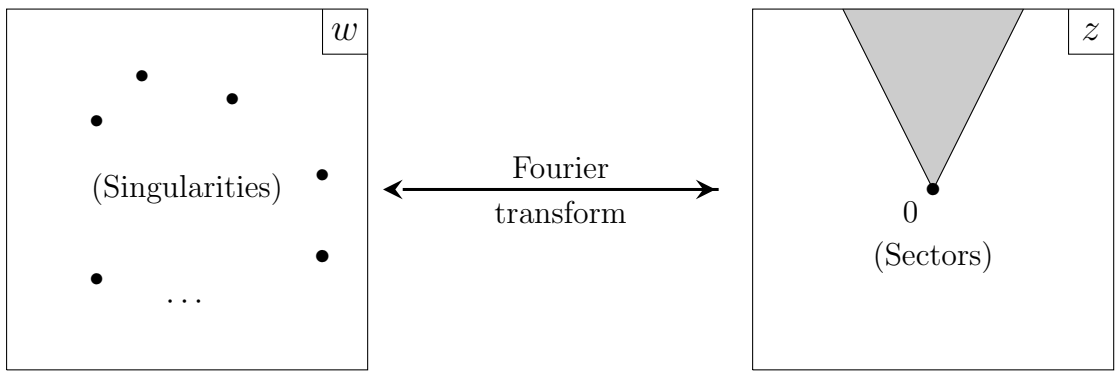

  The formalism of resurgence, see, e.g., \cite{CNP},  is based on the 
 Borel summation process for divergent series which  can be seen as an
 application of the Fourier transform in the complex domain.   It connects
 two copies of the complex plane $\CC$ which are loosely related 
 ``by the Fourier transform'', see Fig. \ref{fig:borel}. 
 
 \vskip .2cm
 
  The main idea is that  complicated ``irregular''  behaviour (sectorial asymptotics
 given by divergent series)
 in the $z$-plane is transformed into ``regular'' behavior (geometry of singularities)
 in the $w$-  (Borel) plane. Note that this is the same idea as studied in
 \S \ref{subsec:fourier}, where we represented some irregular $D$-modules as
 Fourier transforms of regular ones which are given by topological data, and described
 the Stokes matrices in terms of these data. 
 The new phenomenon here is that we typically have infinitely many singularities. 
  This leads to the following philosophy.

  \begin{itemize}
  
 \item[(1)] We need to work in the  plane $\CC = \CC_w$ which is 
  still thought of as the target of the central charge (for the $4$d theory) 
  as well as the Borel plane of the resurgence formalism.   
  
    \item[(2)] We need to study perverse sheaves on this  $\CC$ 
    with possibly  infinitely many singularities.
    
    \item[(3)] The category of such sheaves should have  a symmetric monoidal structure 
     $\star$
    given by  an appropriate version of the (middle) convolution.
    
    \item[(4)] The proper conceptual objects of resurgence theory are perverse sheaves
    which are algebras (commutative or not) with respect to $\star$,
    or modules over such algebras, or possibly  other algebraic structures in $\star$. 
    
    \item[(5)] Note that if $\Fc\in \Perv(\CC,A)$ and $\Gc\in\Perv(\CC,B)$,
    then
    $\Fc\star\Gc\in \Perv(\CC, A+B)$ where $A+B = \{a+b| \, a\in A, b\in B\}$.
    So a for a $\star$-algebra $\Fc$ the singularities typically form a semigroup,
    which must be infinite unless trivial. In examples of interest this semigroup
    is the {\em lattice of periods}, see \cite{Kont-S, Kont-S-Floer}.

    \item[(6)] The Fourier transform changes $\star$ into fiberwise product
    so a $\star$-algebra is transformed into an irregular local system of
    infinite rank with each fiber being an algebra and Stokes transformations
    being homomorphisms of algebras. 
    
    \item[(7)] Finally, the above considerations should eventually be categorified 
    at the level of schobers. 
    
    \end{itemize}

   \noindent There are several technical issues with realization of this program, having to do,
    e.g.,  with
    perverse sheaves possibly having infinite  generic rank (even if all the spaces of vanishing
    cycles are finite-dimensional) as well as with the possibility of the set
    of singularities being a countable but dense subset such as a free abelian subgroup
    in $\CC$ of rank $\geq 3$. They are now subject of work in progress by the authors.

\vfill\eject

	\vskip 1cm

M.K.: Kavli IPMU, 5-1-5 Kashiwanoha, Kashiwa, Chiba, 277-8583 Japan. \hfill\break
 Email: {\tt mikhail.kapranov@protonmail.com}

 \vskip .2cm

 Y.S.: Department of Mathematics, Kansas State University, Manhattan KS 66508 USA.
 \hfill\break Email: {\tt soibel@math.ksu.edu}

\end{document}